\documentclass{amsart}
\usepackage{amssymb} 
\usepackage{amsmath} 
\usepackage{amscd}
\usepackage[matrix,arrow]{xy}

\DeclareMathOperator{\hrm}{h}
\newcommand{\eps}{\varepsilon}

\newcommand{\vv}{{\upsilon}}

\DeclareMathOperator{\norm}{Norm}

\DeclareMathOperator{\ord}{ord}

\newcommand{\Q}{{\mathbb Q}}
\newcommand{\Z}{{\mathbb Z}}

\newcommand{\C}{{\mathbb C}}
\newcommand{\R}{{\mathbb R}}
\newcommand{\F}{{\mathbb F}}

\newcommand{\cK}{\mathcal{K}}

\newcommand{\cB}{\mathcal{B}}

\newcommand{\cE}{\mathcal{E}}

\newcommand{\cM}{\mathcal{M}}

\newcommand{\OO}{{\mathcal O}}

\newcommand{\bl}{\mathbf{l}}
\newcommand{\bx}{\mathbf{x}}

\begin {document}

\newtheorem{thm}{Theorem}
\newtheorem{lem}{Lemma}[section]
\newtheorem{prop}[lem]{Proposition}

\theoremstyle{definition}

\theoremstyle{remark}

\title[]{Integral Points on Hyperelliptic Curves}

\author{Y.\ Bugeaud, M.\ Mignotte, 
S.\ Siksek, M.\ Stoll, Sz.\ Tengely}

\address{Yann Bugeaud and Maurice Mignotte\\
Universit\'{e} Louis Pasteur\\
U. F. R. de math\'{e}matiques\\
7, rue Ren\'{e} Descartes\\
67084 Strasbourg Cedex\\
France}
\email{bugeaud@math.u-strasbg.fr}
\email{mignotte@math.u-strasbg.fr}

\address{Samir Siksek\\
	Institute of Mathematics\\
	University of Warwick\\
	Coventry\\
	CV4 7AL \\
	United Kingdom}

\email{s.siksek@warwick.ac.uk}

\address{Michael Stoll\\
Mathematisches Institut\\ 
Universit\"{a}t Bayreuth\\ 
95440 Bayreuth\\
Germany}

\email{Michael.Stoll@uni-bayreuth.de}

\address{Szabolcs Tengely\\
Institute of Mathematics\\
University of Debrecen
and the Number Theory Research Group
of the Hungarian Academy of Sciences\\
P.O.Box 12\\
4010 Debrecen\\
Hungary}

\email{tengely@math.klte.hu}

\date{\today}
\thanks{}


\begin{abstract}
Let $C: Y^2=a_n X^n+\cdots+a_0$ be a hyperelliptic
curve with the $a_i$ rational integers, $n \geq 5$, and the
polynomial on the right irreducible. Let $J$ be its Jacobian.
We give a completely explicit upper bound for the integral points
on the model $C$, provided we know at least one rational
point on $C$ and a Mordell--Weil basis for $J(\Q)$. 
We also explain a powerful refinement of the Mordell--Weil sieve
which, combined with the upper bound, is capable of determining
all the integral points. Our method is illustrated 
by determining the integral points on the
genus $2$ hyperelliptic models $Y^2-Y=X^5-X$
and $\binom{Y}{2}=\binom{X}{5}$.
\end{abstract}
\maketitle

\section{Introduction}

Consider the hyperelliptic curve with affine model
\begin{equation}\label{eqn:hyper}
C: Y^2=a_n X^n+a_{n-1}X^{n-1}+\cdots+a_0,
\end{equation}
with $a_0,\dotsc,a_n$ rational integers, $a_n \neq 0$, $n \geq 5$,
and the polynomial on the right irreducible.
Let $H=\max \{\lvert a_0 \rvert, \ldots , \lvert a_n \rvert\}$. 
In one of the earliest applications
of his theory of lower bounds for linear forms in logarithms, 
Baker \cite{Baker} showed that 
any integral point $(X,Y)$ on this affine model satisfies
\[
\max(\lvert X\rvert,\lvert Y\rvert) \leq
\exp \exp \exp \{(n^{10n}H)^{n^2} \}.
\]
Such bounds have been improved considerably by many authors,
including 
Sprind\v{z}uk \cite{Sprind}, 
Brindza \cite{Brindza}, Schmidt \cite{Schmidt}, Poulakis \cite{Poulakis},
Bilu \cite{Bilu}, Bugeaud \cite{Bu1} and Voutier \cite{Voutier}.  
Despite the improvements, the bounds remain astronomical and often
involve inexplicit constants. 

In this paper we explain a new 
method for explicitly computing the integral points on affine
models of hyperelliptic curves \eqref{eqn:hyper}. The method falls into
two distinct steps:
\begin{enumerate}
\item[(i)] We give a completely explicit upper bound for the size
of integral solutions of \eqref{eqn:hyper}. This upper bound
combines the many refinements found in the papers of Voutier,
Bugeaud, etc., together with Matveev's  
bounds for linear forms in logarithms \cite{Matveev},
and a method for bounding the regulators based on a theorem
of Landau \cite{Landau}. 
\item[(ii)] The bounds obtained in (i), whilst substantially
better than bounds given by earlier authors, are still astronomical.
We explain a powerful variant of the Mordell--Weil sieve
which, combined with the bound obtained in (i),
 is capable of showing that the known solutions to~\eqref{eqn:hyper}
are the only ones. 
\end{enumerate}
Step (i) requires two assumptions:
\begin{enumerate}
\item[(a)] We assume that we know at least one
rational point $P_0$ on $C$.
\item[(b)] Let $J$ be the Jacobian of $C$. We assume
that a Mordell--Weil basis for $J(\Q)$ is known.
\end{enumerate}
For step (ii) we need assumptions (a), (b) and also:
\begin{enumerate}
\item[(c)] We assume that the canonical height $\hat{h}:J(\Q)\rightarrow \R$
is explicitly computable and that we have explicit bounds
for the difference
\begin{equation}\label{eqn:hdiff}
\mu_1 \leq h(D)-\hat{h}(D) \leq \mu_1^\prime
\end{equation}
where $h$ is an appropriately normalized logarithmic height on $J$
that allows us to enumerate points~$P$ in~$J(\Q)$ with $h(P) \le B$
for a given bound~$B$.
\end{enumerate}
Assumptions (a)--(c) deserve a comment or two. 
For many families of curves of higher genus, practical descent
strategies are available for estimating the
 rank of the Mordell--Weil group; see for example
\cite{CF},
\cite{PS}, \cite{Sch} and
\cite{Stoll}.
To provably determine the Mordell--Weil group 
one however needs bounds for the difference between the
logarithmic and canonical heights. For Jacobians of curves
of genus $2$ such bounds have been determined by Stoll \cite{StollB1},
\cite{StollB2}, building on previous work of Flynn and Smart \cite{FS}.  
At present, no such bounds have been determined for Jacobians
of curves of genus $\geq 3$, although work on this is in progress. 
The assumption about the knowledge of a rational point is
a common sense assumption that brings some simplifications
to our method,
although the method can be modified to cope with the situation
where no rational point is known. However, if a search on
a curve of genus $\geq 2$ 
reveals no rational points, 
it is probable that there are none, and the methods
of \cite{BS},  \cite{BS2}, \cite{BS3} are  
likely to succeed in  proving this.

We illustrate the practicality of our approach by proving the following
results.
\begin{thm}\label{thm:ex}
The only integral solutions to the equation
\begin{equation}\label{eqn:ex}
Y^2-Y=X^5-X
\end{equation}
are 
\begin{align*}
(X,Y) &= (-1,0),\; (-1,1), \; (0,0),\; (0,1),\; (1,0),\; (1,1),\;(2,-5),\\
     & \qquad (2,6),\; (3,-15),\; (3,16),\; (30,-4929),\; (30,4930).
\end{align*}
\end{thm}
\begin{thm}\label{thm:ex2}
The only integral solutions to the equation
\begin{equation}\label{eqn:ex2}
\binom{Y}{2}=\binom{X}{5}
\end{equation}
are
\begin{gather*}
(X,Y)=(0, 0), \;
(0, 1),\;
(1, 0),\;
(1, 1),\;
(2, 0),\;
(2, 1),\;
(3, 0),\;
(3, 1),\;
(4, 0),\;
(4, 1),\;
(5, -1),\\
(5, 2),\;
(6, -3),\;
(6, 4),\;
(7, -6),\;
(7, 7),\;
(15, -77),\;
(15, 78),\;
(19, -152),\;
(19, 153).
\end{gather*}
\end{thm}
Equations \eqref{eqn:ex} and \eqref{eqn:ex2} are
of historical interest and Section~\ref{sec:history}
gives a brief outline of their history. 
For now we merely mention that these two equations
are  the first two problems on a list
of 22 unsolved Diophantine problems \cite{ET}, compiled
by Evertse and Tijdeman following a recent workshop on Diophantine
equations at Leiden.

To appreciate why the innocent-looking equations \eqref{eqn:ex}
and \eqref{eqn:ex2}
have resisted previous attempts, let 
us briefly
survey the available methods which apply to hyperelliptic curves
and then briefly explain why they fail in these cases. To determine
the integral points on the affine model $C$ given by 
an equation \eqref{eqn:hyper} there are four available methods: 
\begin{enumerate}
\item[(I)] The first is Chabauty's elegant method which in fact determines
all rational points on $C$ 
in many cases,
provided the rank of the Mordell--Weil
group of its Jacobian is strictly less than the genus $g$;
see for example \cite{Flynnc}, 
\cite{We}.
Chabauty's method fails if the rank of the Mordell--Weil
group exceeds the genus.

\item[(II)] A second method is to use coverings, often combined with a version
of Chabauty called 
\lq Elliptic Curve Chabauty\rq. 
See 
\cite{Br1}, \cite{Br2}, \cite{FW1}, \cite{FW2}.
 This approach often requires computations of Mordell--Weil groups
 over number fields (and does fail if the rank of the Mordell--Weil
 groups is too large). 

\item[(III)] A third method is to combine Baker's approach through 
$S$-units with the LLL
algorithm to obtain all the solutions provided that certain 
relevant unit groups and class groups can be computed; for
a modern treatment, see \cite{BH} or \cite[Section XIV.4]{Smart}. 
This strategy often fails
in practice as the number fields involved have very high degree.

\item[(IV)] A fourth 
approach
is to apply Skolem's method to the $S$-unit equations
(see \cite[Section III.2]{Smart}).
This needs the same expensive information as the third method.
\end{enumerate}
The Jacobians of the curves given by \eqref{eqn:ex} and \eqref{eqn:ex2}
respectively have ranks $3$ and $6$
and so Chabauty's method fails.
To employ 
Elliptic Curve Chabauty 
would require the computation of 
Mordell--Weil groups of elliptic curves 
without rational 2-torsion
over number fields of degree $5$
(which does not seem practical at present).
To apply the $S$-unit approach (with either LLL or Skolem)
requires the computations of the unit groups and class groups of 
several number fields of degree $40$; 
a computation that seems 
completely impractical at present.

\medskip

Our paper is arranged as follows. Section \ref{sec:history}
gives a brief history of equations~\eqref{eqn:ex}
and \eqref{eqn:ex2}.
In Section \ref{sec:descent} we show, 
after appropriate scaling, that 
an integral point $(x,y)$ satisfies
$x-\alpha=\kappa \xi^2$ where $\alpha$
is some fixed algebraic integer,
$\xi \in \Q(\alpha)$,
and $\kappa$ is an algebraic integer
belonging to a finite computable set. In Section \ref{sec:bounds}
we give bounds for the size of solutions $x \in \Z$ to an equation
of the form $x-\alpha=\kappa \xi^2$ where $\alpha$ and $\kappa$
are fixed algebraic integers. Thus, in effect, we obtain bounds for
the size of solutions integral points on our affine model
for \eqref{eqn:hyper}. Sections \ref{sec:heights}--\ref{sec:uniteq}
are preparation for Section \ref{sec:bounds}: in particular Section \ref{sec:heights}
is concerned with heights; Section \ref{sec:regbounds} explains how a 
theorem 
of Landau can be used to bound the regulators of number fields;
Section \ref{sec:fu} collects and refines various results on appropriate
choices of systems of fundamental units; Section \ref{sec:Matveev}
is devoted to Matveev's bounds for linear forms in logarithms;
in Section \ref{sec:uniteq} we use Matveev's bounds and the results
of previous sections to prove a bound on the size of solutions of unit
equations; in Section \ref{sec:bounds} we deduce the bounds
for $x$ alluded to above from the bounds for solutions of unit equations.
Despite our best efforts, the bounds obtained for $x$ are still so large
that no naive search up to those bounds is conceivable. Over the next
three sections \ref{sec:mwI}, \ref{sec:mwII}, \ref{sec:lb}
we explain how to sieve effectively up to these bounds using the
Mordell--Weil group of the Jacobian. In particular, Section~\ref{sec:mwII}
gives a powerful refinement of the 
Mordell--Weil sieve
(\cite{BS},  \cite{BS3})
which we expect to have applications elsewhere. Finally,
in Section~\ref{sec:example} we apply the method of this paper
to prove Theorems~\ref{thm:ex} and~\ref{thm:ex2}.
 
We are grateful to the referee and editors for many useful comments,
and to Mr.\ Homero Gallegos--Ruiz for spotting many misprints. 

\section{History of Equations~\eqref{eqn:ex} and \eqref{eqn:ex2}}
\label{sec:history}
The equation \eqref{eqn:ex} is a special case of the 
family of Diophantine
equations
\begin{equation}\label{eqn:FA}
Y^p-Y=X^q-X, \qquad 2 \leq p < q.  
\end{equation}
This family has previously been studied by Fielder and Alford \cite{FA}
and by Mignotte and Peth\H{o} \cite{MigP}. 
The (genus $1$) case $p=2$, $q=3$
was solved by Mordell \cite{Mordell} who showed that the only solutions
in this case are
\begin{gather*}
(X,Y)=(0,0),(0,1),(\pm 1,0),(\pm 1,1),(2,3), (2,-2), (6,15),(6,-14).
\end{gather*}
Fielder and Alford  presented the following 
list of solutions with $X$, $Y>1$:
\begin{eqnarray*}
(p,q,X,Y)=(2, 3, 2, 3),
(2, 3, 6, 15),
(2, 5, 2, 6),
(2, 5, 3, 16),\\
(2, 5, 30, 4930),
(2, 7, 5, 280),
(2, 13, 2, 91),
(3, 7, 3, 13).
\end{eqnarray*}
Mignotte and Peth\H{o} proved that for 
given $p$ and $q$ with $2\leq p<q$, the Diophantine 
equation \eqref{eqn:FA} has only a finite number of integral solutions. 
Assuming the $abc$-conjecture, they showed that 
equation \eqref{eqn:FA} has only finitely many 
solutions with $X$, $Y>1$. 

If $p=2$, $q>2$ and $y$ is a prime power, 
then Mignotte and Peth\H{o} found all solutions of the equation
and these are all in  Fielder and Alford's list.

\medskip

Equation~\eqref{eqn:ex2} is a special case
of the Diophantine equation 
\begin{equation}\label{binom}
\binom{n}{k}=\binom{m}{l},
\end{equation}
in unknowns $k$, $l$, $m$, $n$.
This is usually considered with the restrictions
$2\leq k\leq n/2$, and $2\leq l\leq m/2$.
The only known  solutions (with these restrictions)
are the following
\begin{eqnarray*}
&& \binom{16}{2}=\binom{10}{3},\quad \binom{56}{2}=\binom{22}{3},\quad \binom{120}{2}=\binom{36}{3},\\
&& \binom{21}{2}=\binom{10}{4},\quad \binom{153}{2}=\binom{19}{5},\quad \binom{78}{2}=\binom{15}{5}=\binom{14}{6},\\
&& \binom{221}{2}=\binom{17}{8},\quad \binom{F_{2i+2}F_{2i+3}}{F_{2i}F_{2i+3}}=\binom{F_{2i+2}F_{2i+3}-1}{F_{2i}F_{2i+3}+1} \mbox{ for } i=1,2,\ldots,
\end{eqnarray*}
where $F_n$ is the $n$th Fibonacci number.
It is known that there are no other non-trivial 
solutions with $\binom{n}{k}\leq 10^{30}$ or $n\leq 1000$;
see \cite{deWeger2}. 
The infinite family of solutions was found by 
Lind \cite{Lind} and Singmaster \cite{Singmaster}.

Equation \eqref{binom} has been completely solved for
pairs 
\[
(k,l)=(2,3),\; (2,4),\; (2,6),\; (2,8),\; (3,4),\; (3,6),\; (4,6), \; (4,8).
\] 
These are the cases when one can easily reduce 
the equation to the determination of solutions of a number 
of Thue equations or elliptic Diophantine equations.
In 1966, Avanesov \cite{Avanesov} found all solutions of 
equation \eqref{binom} with $(k,l)=(2,3)$. 
De~Weger \cite{deWeger} and independently 
Pint\'er \cite{Pinter} solved the equation with $(k,l)=(2,4)$. 
The case $(k,l)=(3,4)$ reduces to the equation 
$Y(Y+1)=X(X+1)(X+2)$ which was solved by Mordell \cite{Mordell}.
The remaining pairs $(2,6),(2,8),(3,6),(4,6), (4,8)$ were 
treated by Stroeker and de Weger \cite{StroekerWeger}, 
using linear forms in elliptic logarithms.

There are also some general finiteness results 
related to equation \eqref{binom}. 
In 1988, Kiss \cite{Kiss} 
proved that if $k=2$ and $l$ is a 
given odd prime, then the equation has only finitely 
many positive integral solutions. 
Using Baker's method, Brindza \cite{Brindza2} 
showed that equation \eqref{binom} with $k=2$ and $l\geq 3$ has 
only finitely many positive integral solutions.

\section{Descent}\label{sec:descent}
Consider the integral points on the affine model
of the hyperelliptic curve \eqref{eqn:hyper}. 
If the polynomial on the right-hand side is reducible then 
the obvious factorisation
argument reduces the problem of determining the integral points
on \eqref{eqn:hyper} to determining those on simpler hyperelliptic
curves, or on genus $1$ curves. The integral
points on a genus $1$ curve can be determined by highly
successful algorithms
(e.g.\ 
\cite{Smart},
\cite{sttz2})
based on LLL and David's bound  for linear forms 
in elliptic logarithms. 

We therefore suppose henceforth that the 
polynomial on the right-hand side of \eqref{eqn:hyper}
is irreducible; this is certainly the most difficult case.
By appropriate scaling, one transforms the problem of integral
points on \eqref{eqn:hyper} to integral points on a model
of the form
\begin{equation}\label{eqn:twist}
a y^2=x^n+b_{n-1}x^{n-1}+\dots+b_0,
\end{equation}
where $a$ and the $b_i$ are integers, with $a \neq 0$. 
We shall work henceforth with this model of the hyperelliptic
curve.
Denote the polynomial on the right-hand side by $f$
and let $\alpha$ be a root of $f$.
Then a standard argument shows that
\[
x-\alpha=\kappa \xi^2
\]
where $\kappa$, $\xi \in K=\Q(\alpha)$
and $\kappa$ is an {\bf algebraic integer that comes from 
a finite computable set}. In this section we suppose 
that the Mordell--Weil group $J(\Q)$ of the curve $C$
is known, and we show how to compute such a set
of $\kappa$ using our knowledge of the Mordell--Weil
group $J(\Q)$. The method for doing this depends on whether
the degree $n$ is odd or even. 

\subsection{The Odd Degree Case}

Each coset of $J(\Q)/2J(\Q)$ has a coset representative
of the form $\sum_{i=1}^{m} (P_i -\infty)$ where the set 
$\{P_1,\dotsc,P_m\}$
is stable 
under the action of Galois, and
 where all $y(P_i)$ are non-zero. 
Now write $x(P_i)=\gamma_i/d_i^2$ where $\gamma_i$
is an algebraic integer and $d_i \in \Z_{\geq 1}$; moreover if 
$P_i$, $P_j$ are conjugate then we may suppose that $d_i=d_j$
and so $\gamma_i$, $\gamma_j$ are conjugate. To such a coset
representative of $J(\Q)/2J(\Q)$ we associate 
\[
\kappa=a^{(\text{$m$ mod $2$})}\prod_{i=1}^{m}\left(\gamma_i-\alpha d_i^2 \right).
\]
\begin{lem}\label{lem:cK}
Let $\cK$ be a set of  $\kappa$ associated as above to a complete 
set of coset representatives of $J(\Q)/2J(\Q)$. 
Then $\cK$ is a finite subset of  $\OO_K$ 
and if $(x,y)$ is
an integral point on the model \eqref{eqn:twist}
then $x-\alpha=\kappa \xi^2$
for some $\kappa \in \cK$ and $\xi \in K$. 
\end{lem}
\begin{proof}
This follows trivially from the standard homomorphism
\[
\theta : J(\Q)/2J(\Q) \rightarrow K^*/{K^*}^2
\]
that is given by 
\[
\theta\left(\sum_{i=1}^m (P_i-\infty)\right)=
a^m \prod_{i=1}^m \left(x(P_i)-\alpha \right) \pmod{{K^*}^2}
\]
for coset representatives $\sum (P_i -\infty)$ with $y(P_i) \neq 0$;
see Section 4 of \cite{Stoll}.
\end{proof}

\subsection{The Even Degree Case} 
As mentioned in the introduction, we shall 
assume 
the existence of at least one rational point $P_0$.
If $P_0$ is one of the two points at infinity, let 
$\epsilon_0=1$. Otherwise, 
as $f$ is irreducible, $y(P_0) \neq 0$; write $x(P_0)=\gamma_0/d_0^2$
with $\gamma_0 \in \Z$ and $d_0 \in \Z_{\geq 1}$
and let $\epsilon_0=\gamma_0-\alpha d_0^2$. 

Each coset of $J(\Q)/2J(\Q)$ has a coset representative
of the form $\sum_{i=1}^{m} (P_i -P_0)$ where the set 
$\{P_1,\dotsc,P_m\}$
is stable 
under the action of Galois, and
 where all $y(P_i)$ are non-zero for $i=1,\dotsc,m$. 
Write $x(P_i)=\gamma_i/d_i^2$ where $\gamma_i$
is an algebraic integer and $d_i \in \Z_{\geq 1}$; moreover if 
$P_i$, $P_j$ are conjugate then we may suppose that $d_i=d_j$
and so $\gamma_i$, $\gamma_j$ are conjugate. To such a coset
representative of $J(\Q)/2J(\Q)$ we associate 
\[
\epsilon=\epsilon_0^{(\text{$m$ mod $2$})}\prod_{i=1}^{m}\left(\gamma_i-\alpha d_i^2 \right).
\]
\begin{lem}\label{lem:cK2}
Let $\cE$ be a set of  $\epsilon$ associated as above to a complete 
set of coset representatives of $J(\Q)/2J(\Q)$.
Let $\Delta$ be the discriminant of the polynomial $f$.
For each $\epsilon \in \cE$,
let $\cB_\epsilon$ be the set of square-free rational integers 
supported only by primes dividing $a \Delta \norm_{K/\Q}(\epsilon)$.  
Let $\cK=\{ \epsilon b : \epsilon \in \cE, \; b \in \cB_\epsilon\}$.
Then $\cK$ is a finite subset of $\OO_K$ and if $(x,y)$ is
an integral point on the model \eqref{eqn:twist}
then $x-\alpha=\kappa \xi^2$
for some $\kappa \in \cK$ and $\xi \in K$. 
\end{lem}
\begin{proof}
In our even degree case, 
the homomorphism~$\theta$ takes values in $K^*/\Q^* {K^*}^2$.
Thus if $(x,y)$ is an integral point on the model~\eqref{eqn:twist},
we have that $(x-\alpha)=\epsilon b \xi^2$ for some $\epsilon \in \cE$
and $b$ a square-free rational integer. 
A standard argument
shows that $2 \mid \ord_\wp(x-\alpha)$ for all prime ideals 
$\wp \nmid a\Delta$. Hence,  $2 \mid \ord_\wp(b)$
for all $\wp \nmid a \Delta \epsilon$. Let $\wp \mid p$
where $p$ is a rational prime not dividing 
$a \Delta \norm_{K/\Q}(\epsilon)$. Then $p$ is unramified
in $K/\Q$ and so $\ord_p(b)=\ord_\wp(b) \equiv 0 \pmod{2}$.
This shows that $b \in \cB_\epsilon$ and proves the lemma.
\end{proof}


\subsection{Remarks} The following remarks are applicable
both to the odd and the even degree cases. 
\begin{itemize}
\item We point out that even if we do not know coset representatives
for $J(\Q)/2J(\Q)$, we can still obtain a suitable 
(though larger) set of $\kappa$
that satisfies the conclusions of Lemmas~\ref{lem:cK} 
and~\ref{lem:cK2} provided we are
able to compute
the class group and unit group of the number field $K$; for
this see for example \cite[Section 2.2]{Br1}.
\item We can use local information
at small and bad primes to restrict the set $\cK$ further, 
compare~\cite{BS}
and~\cite{BS2}, where this is applied to rational points. In our case,
we can restrict the local computations to $x \in \Z_p$ instead of~$\Q_p$.
\end{itemize}

\section{Heights}\label{sec:heights}

We fix once and for all the following notation.

\medskip
\begin{tabular}{lll}
$K$ & $\qquad$ & {a number field,}\\
$\OO_K$ &  & {the ring of integers of $K$,}\\
$M_K$ &  & {the set of all places of $K$,} \\
$M_K^0$ & & {the set of non-Archimedean places of $K$,}\\
$M_K^\infty$ & & {the set of Archimedean places of $K$,}\\
$\vv$ & & {a place of $K$,}\\
$K_\vv$ & & {the completion of~$K$ at~$\vv$,}\\
$d_\vv$ & & {the local degree $[K_\vv:\Q_\vv]$.}
\end{tabular}
\medskip

For $\vv \in M_K$, we let $\lvert \cdot \rvert_\vv$ be 
the usual normalized
valuation corresponding to $\vv$; in particular
if $\vv$ is non-Archimedean and $p$ is the rational
prime below $\vv$ then $\lvert p \rvert_\vv=p^{-1}$. 
Thus if $L/K$ is a field extension, 
and $\omega$ a place of $L$ above $\vv$ then 
$\lvert \alpha\rvert_\omega=\lvert \alpha \rvert_\vv$,
for all $\alpha \in K$.

Define 
\[
\lVert \alpha \rVert_\vv = \lvert \alpha \rvert_\vv^{d_\vv}.
\]
Hence for $\alpha \in K^*$, the product formula states that
\[
\prod_{\vv \in M_K} \lVert \alpha \rVert_\vv=1.
\]
In particular, if $\vv$ is Archimedean, corresponding to a real or complex
embedding $\sigma$ of $K$ then 
\[
\lvert \alpha \rvert_\vv=\lvert \sigma(\alpha) \rvert 
\qquad \text{and} \qquad
\lVert \alpha \rVert_\vv=
\begin{cases}
\lvert \sigma(\alpha) \rvert & \text{if $\sigma$ is real} \\
\lvert \sigma(\alpha) \rvert^2 & \text{if $\sigma$ is complex}.
\end{cases}
\]

For $\alpha \in K$, the (absolute) logarithmic height $\hrm(\alpha)$
is given by
\begin{equation}\label{eqn:loght}
        \hrm(\alpha)= 
	\frac{1}{[K:\Q]} \sum_{\vv \in M_K } d_\vv \log \max \left\{ 1 ,
        \lvert \alpha \rvert_\vv \right\}=
	\frac{1}{[K:\Q]} \sum_{\vv \in M_K } \log \max \left\{ 1 ,
        \lVert \alpha \rVert_\vv \right\}.
\end{equation}
The absolute logarithmic height of $\alpha$ 
is independent of the field $K$ containing $\alpha$. 

We shall need the following elementary properties of heights.
\begin{lem}\label{lem:elementary}
For any non-zero algebraic number $\alpha$, we have 
$\hrm(\alpha^{-1})=\hrm(\alpha)$.
For algebraic numbers $\alpha_1,\dotsc,\alpha_n$, we have 
\[
\hrm(\alpha_1\alpha_2\cdots\alpha_n) \leq \hrm(\alpha_1)+\dots +\hrm(\alpha_n),
\qquad
\hrm(\alpha_1+\dots+\alpha_n) \leq \log{n}+\hrm(\alpha_1)+\dots +\hrm(\alpha_n).
\]
\end{lem}
\begin{proof}
The lemma is Exercise 8.8 {} in \cite{SilvI}. We do not know
of a reference for the proof and so we will indicate briefly
the proof of the second (more difficult) inequality.
For $\vv \in M_K$, choose $i_\vv$ in $\{1, \ldots , n\}$ to satisfy
$\max\{\lvert \alpha_1 \rvert_\vv,\dotsc, 
\lvert \alpha_n \rvert_\vv \}=\lvert \alpha_{i_\vv} \rvert_\vv$.
Note that
\[
\lvert \alpha_1+\dots+\alpha_n \rvert_\vv
\leq \epsilon_\vv \lvert \alpha_{i_\vv} \rvert_\vv,
\qquad \text{where} \qquad
\epsilon_\vv=\begin{cases}
n & \text{if $\vv$ is Archimedean,} \\
1 & \text{otherwise.}
\end{cases}
\]
Thus
\[
\log \max\{1, \lvert \alpha_1+\dots+\alpha_n \rvert_\vv \}
\leq \log{\epsilon_\vv}+ \log \max\{1, \lvert \alpha_{i_\vv} \rvert_\vv\}
\leq \log{\epsilon_\vv}+ \sum_{i=1}^n \log \max\{1, \lvert \alpha_i \rvert_\vv\}.  
\]
Observe that 
\[
\frac{1}{[K:\Q]}\sum_{\vv \in M_K} d_\vv \log{\epsilon_\vv}=
\frac{\log{n}}{[K:\Q]} \sum_{\vv \in M_K^\infty} d_\vv=
\log{n};
\]
the desired inequality follows from the definition of logarithmic
height \eqref{eqn:loght}.
\end{proof}
\bigskip

\subsection{Height Lower Bound}
We need the following result of 
Voutier \cite{Voutier2} concerning Lehmer's problem.

\begin{lem}\label{lem:Mahler}
Let $K$ be a number field of degree $d$.
Let
\[
\partial_K=\begin{cases}
\frac{\log{2}}{d} & \text{if $d=1$, $2$}, \\
\frac{1}{4}\,\left(\frac{\log \log d}{\log d}\right)^3 & \text{if $d \geq 3$}.
\end{cases}
\]
Then, for every 
non-zero algebraic number $\alpha$ in $K$, which is not
a root of unity, 
\[
\deg(\alpha) \,\hrm (\alpha)\geq \partial_K.
\]
\end{lem}

Throughout, by the logarithm of a complex number, we mean 
the principal determination of the logarithm. In other words,
if $x \in \C^*$ we express $x=r e^{i \theta}$ where
$r >0$ and $-\pi < \theta \leq \pi$;  we then let
$\log{x}=\log{r}+i \theta$.
\begin{lem}\label{lem:abslog}
Let $K$ be a number field and let
\[
\partial^\prime_K=\left( 1 + \frac{\pi^2}{\partial_K^2} \right)^{1/2}.
\]
For any non-zero $\alpha \in K$ 
and any place $\vv \in M_K$
\[
\log \lvert \alpha \rvert_\vv \leq \deg(\alpha) \hrm(\alpha), \qquad
\log \lVert \alpha \rVert_\vv \leq [K:\Q] \hrm(\alpha).
\]
Moreover, if $\alpha$ is not a root of unity and $\sigma$
is a real or complex embedding of $K$ then
\[
\lvert \log \sigma (\alpha) \rvert \leq  \partial^\prime_K \deg(\alpha) \hrm(\alpha).
\]
\end{lem}
\begin{proof}
The first two inequalities are an immediate consequence of the definition
of absolute logarithmic height. 
For the last, write $\sigma(\alpha)  = e^{a + ib}$, with $a=\log \lvert \sigma(\alpha)\rvert$
and $\lvert b\rvert \le \pi$, and let $d=\deg(\alpha)$. 
Then we have
\[
\lvert \log \sigma(\alpha)\rvert 
=(a^2 + b^2 )^{1/2} \leq
(\log^2 \lvert \sigma(\alpha) \rvert + \pi^2)^{1/2}
\leq ((d \hrm(\alpha ))^2 + \pi^2)^{1/2}.
\]
By Lemma~\ref{lem:Mahler}
we have $d \hrm(\alpha ) \geq \partial_K$, so
\[
\lvert \log \sigma(\alpha) \rvert  
\le d \, \hrm(\alpha ) \,
\left( 1 + \frac{\pi^2}{\partial_K^2} \right)^{1/2},
\]
as required.
\end{proof}

\section{Bounds for Regulators} \label{sec:regbounds}

Later on we need to give upper bounds for the regulators of 
complicated number fields of high degree. The following
lemma, based on bounds of Landau \cite{Landau}, is an easy way to obtain
reasonable bounds.

\begin{lem}\label{lem:Landau}
Let $K$ be a number field
with degree $d=u+2v$
where $u$ and $v$ are respectively the numbers of real and complex embeddings.
Denote the absolute discriminant by $D_K$ and the regulator by $R_K$,
and the number of roots of unity in $K$ by $w$. Suppose, moreover,
that $L$ is a real number such that $D_K \leq L$.
Let
\[
            a = 2^{-v} \, \pi^{-d/2} \, \sqrt{L}.
\]
Define the function $f_{K}(L,s)$ by
\[
          f_{K}(L,s)=2^{-u} \, w \,
          a^s \, \bigl( \Gamma (s/2) \bigr)^{u} \, \bigl( \Gamma(s) \bigr)^{v}
          s^{d+1} \, (s-1)^{1-d},
\]
and let $B_{K}(L)=\mathrm{min} \left\{f_{K}(L,2-t/1000)~: t=0,1,\ldots,999 \right\}$.
Then $R_{K} < B_K(L)$.
\end{lem}
\begin{proof}
Landau \cite[proof of Hilfssatz 1]{Landau} established the inequality
$R_{K}< f_K(D_K,s)$
for all $s>1$. It is thus
clear that $R_{K} < B_{K}(L)$.
\end{proof}

Perhaps a comment is in order. For a complicated number
field of high degree it is difficult to calculate the discriminant
$D_{K}$ exactly, though it is easy to give an upper bound $L$ for
its size. It is also difficult to
minimise the function $f_{K}(L,s)$ analytically, but we have
found that the above gives an accurate enough result, which is easy to
calculate on a computer.

\section{Fundamental Units}\label{sec:fu}

For the number fields we are concerned with, we shall need to work
with a certain system of fundamental units, given by the following
lemma due to Bugeaud and Gy\H{o}ry, which is Lemma 1 of \cite{BG1}.

\begin{lem}\label{lem:fugood}
Let $K$ be a number field of degree $d$ and let $r=r_K$ be its unit rank
and $R_K$ its regulator.
Define the constants
\[
c_1 = c_1(K) =\frac{(r\,!)^2}{2^{r-1} d^r} \, , \qquad
c_2= c_2(K) =c_1 \left( \frac{d}{\partial_K}\right)^{r-1} , \qquad
c_3 = c_3(K) =c_1 \frac{d^r}{\partial_K} \, .
\]
Then $K$ admits a system
$\{\eps_1, \dotsc, \eps_r\}$ of fundamental units such that:
\begin{enumerate}
\item[(i)]
$\displaystyle \qquad\qquad\qquad \prod_{i=1}^r \hrm (\eps_i) \leq c_1 R_K$, 

\medskip

\item[(ii)]
$\displaystyle \qquad\qquad\qquad \hrm (\eps_i) \leq c_2 R_K, \quad 1\leq i\leq r$,

\medskip

\item[(iii)]
Write $\cM$ for the $r \times r$-matrix $(\log \lVert \eps_i\rVert_\vv)$
where $\vv$ runs over $r$ of the Archimedean places of $K$ and 
$1 \leq i \leq r$. Then the absolute values of the entries
of $\cM^{-1}$ are bounded above by $c_3$.
\end{enumerate}
\end{lem}

\begin{lem}\label{lem:Best}
Let $K$ be a number field of degree $d$, and  let $\{\eps_1,\dotsc,\eps_r\}$
be a system of fundamental units as in Lemma~\ref{lem:fugood}. Define
the constant 
$c_4=c_4(K)=r d c_3$.
Suppose $\eps=\zeta \eps_1^{b_1} \dots \eps_r^{b_r}$,
where $\zeta$ is a root of unity in $K$. Then 
\[
\max\{\lvert b_1\rvert,\dots,\lvert b_r\rvert \} \leq c_4 \hrm(\eps).
\]
\end{lem}
\begin{proof}
Note that for any Archimedean place $v$ of $K$, 
\[
\log\lVert \eps \rVert_v=\sum b_i \log \lVert \eps_i \rVert_v.
\]
The lemma now follows from part (iii) of Lemma~\ref{lem:fugood}, plus
the fact that $\log \lVert \eps \rVert_v \leq d \hrm(\eps)$ for all $v$
given by Lemma~\ref{lem:abslog}.
\end{proof}

\medskip

The following result is a special case of Lemma 2 of \cite{BG1}.

\begin{lem}\label{lem:improve}
Let $K$ be a number field of unit rank $r$ and regulator $K$.
Let $\alpha$ be a non-zero algebraic integer 
belonging to $K$. 
Then there exists a  unit  
$\eps$ of $K$ such that
\[
\hrm (\alpha \eps)\leq c_5 R_K+\frac{\log \lvert \norm_{K/\Q}(\alpha)\rvert}{[K:\Q]}
\]
where
\[
c_5 = c_5(K)= \frac{r^{r+1}}{2  \partial_K^{r-1}}.
\]
\end{lem}

\begin{lem}\label{lem:minsig}
Let $K$ be a number field, $\beta$, $\eps \in K^*$ with $\eps$
being a unit. Let $\sigma$ be the real or complex embedding that
makes $\lvert \sigma(\beta \eps)\rvert$ minimal. Then
\[
\hrm (\beta \eps) \leq \hrm(\beta) - \log \lvert \sigma(\beta \eps)\rvert.
\]
\end{lem}
\begin{proof}
As usual, write $d=[K:\Q]$ and $d_\vv=[K_\vv:\Q_\vv]$.
Note
\begin{equation*}
\begin{split}
\hrm(\beta \eps) & = \hrm(1/\beta\eps) \\
& =\frac{1}{d} \sum_{\vv \in M_K^\infty} d_\vv \max \{0, \log (\lvert \beta \eps \rvert_\vv^{-1})\}
+ \frac{1}{d}\sum_{\vv \in M_K^0} d_\vv \max \{0, \log (\lvert \beta \eps \rvert_\vv^{-1})\} \\
& \leq  \log (\lvert \sigma(\beta \eps)\rvert^{-1})+\frac{1}{d} \sum_{\vv \in M_K^0} d_\vv \max \{0, \log (\lvert \beta \rvert_\vv^{-1})\} \\
&\leq  - \log \lvert \sigma(\beta \eps) \rvert + \frac{1}{d}\sum_{\vv \in M_K} d_\vv \max \{0, \log (\lvert \beta \rvert_\vv^{-1})\} \\
&\leq  - \log \lvert \sigma(\beta \eps) \rvert + \hrm(\beta), 
\end{split}
\end{equation*}
as required.
\end{proof}

\section{Matveev's Lower Bound for Linear Forms in Logarithms}\label{sec:Matveev}

Let $L$ be a number field and let $\sigma$ be a real or complex embedding.
For $\alpha \in L^*$ we define the 
{\em modified
logarithmic height of $\alpha$ with respect to $\sigma$} to be
\[
\hrm_{L,\sigma}(\alpha):=\max \{[L:\Q] \hrm(\alpha) \; , \; \lvert \log \sigma(\alpha) \rvert
\; ,\; 0{.}16 \}.
\]
The modified height is clearly dependent on the number field; we shall need the 
following Lemma which gives a relation between the modified and absolute
height.

\begin{lem}\label{lem:modified}
Let $K \subseteq L$ be number fields and write
\[
\partial_{L/K}= \max \left\{[L:\Q]\; , \; [K:\Q] \partial^\prime_K
\; , \; \frac{0.16[K:\Q]}{\partial_K} \right\}.
\]
Then for any $\alpha \in K$ which is neither zero nor a root of unity,
and any real or complex embedding $\sigma$ of $L$,
\[
\hrm_{L,\sigma}(\alpha) \leq \partial_{L/K} \hrm(\alpha).
\]
\end{lem}
\begin{proof}
By Lemma~\ref{lem:abslog} we have 
\[
[K:\Q] \partial^\prime_K \hrm(\alpha) \geq \partial^\prime_K \deg(\alpha) \hrm(\alpha)
\geq \lvert \log \sigma(\alpha) \rvert.
\]
Moreover, by Lemma~\ref{lem:Mahler},
\[
\frac{0{.}16[K:\Q] \hrm(\alpha)}{\partial_K} \geq 
\frac{0{.}16 \deg(\alpha) \hrm(\alpha)}{\partial_K} \geq 0{.}16.
\]
The lemma follows.
\end{proof} 
\medskip

We shall apply lower bounds on linear forms, more
precisely a version of Matveev's   
estimates \cite{Matveev}.
We recall that 
$\log$ denotes the principal
determination of the logarithm.

\begin{lem}\label{lem:Matveev}
Let $L$ be a number field of degree $d$, 
with $\alpha_1,\dotsc,\alpha_n\in L^*$.
Define 
a
constant 
\[
C(L,n):= 3 \cdot 30^{n+4}\cdot (n+1)^{5{.}5}\,d^2 \,(1+\log{d}).
\]
Consider the \lq\lq linear form\rq\rq\
\[
\Lambda:= \alpha_1^{b_1}\cdots \alpha_n^{b_n}-1,
\]
where $b_1,\dotsc,b_n$ are rational integers and let 
$B:=\max\{\lvert b_1\rvert,\dotsc,\lvert b_n\rvert\}$.
If $\Lambda\neq 0$, and $\sigma$ is any real or complex embedding of $L$ then
\[
\log \lvert \sigma(\Lambda) \rvert > 
-C(L,n) (1+\log(nB)) \prod_{j=1}^n \hrm_{L,\sigma}(\alpha_j).
\]
\end{lem}
\begin{proof}
This 
straightforward corollary of Matveev's estimates is Theorem 9.4 of \cite{BMS1}.
\end{proof}

\medskip

\section{Bounds for Unit Equations}\label{sec:uniteq}
Now we are ready to prove an explicit version of Lemma 4 of \cite{Bu1}.
The proposition below allows us to replace in the final estimate
the regulator of the larger field by the product of the
regulators of two of its subfields. This often results in a
significant improvement of the upper bound for the height.
This idea is due to Voutier \cite{Voutier}.  

\medskip

\begin{prop}\label{prop}
Let 
$L$ 
be a number field of degree $d$, which contains $K_1$ and  
$K_2$ as subfields. 
Let $R_{K_i}$ (respectively $r_i$) be the regulator 
(respectively the unit rank) of $K_i$. 
Suppose further that
$\nu_1$, $\nu_2$ and $\nu_3$ are non-zero elements of $L$ 
with height  
${}\le H$ (with $H\geq 1$) and consider
the unit equation
\begin{equation}\label{eqn:uniteq}
\nu_1\eps_1+\nu_2\eps_2+\nu_3\eps_3 = 0
\end{equation}
where $\eps_1$ is a unit of $K_1$, $\eps_2$ a unit of $K_2$ and $\eps_3$  
a unit of $L$. 
 Then, for $i=1$ and~$2$,
\[
\hrm(\nu_i\eps_i/\nu_3\eps_3) \leq A_2+A_1 
\log \{H+\max\{\hrm(\nu_1 \eps_1),\hrm(\nu_2 \eps_2)\}\},
\]
where 
\[
A_1=2 H \cdot C(L,r_1+r_2+1) 
\cdot c_1(K_1) c_1(K_2) \partial_{L/L}
\cdot
(\partial_{L/K_1})^{r_1}
\cdot
(\partial_{L/K_2})^{r_2}
\cdot R_{K_1} R_{K_2},
\] 
and
\[
A_2=2H+A_1+A_1 \log\{(r_1+r_2+1) \cdot \max\{c_4(K_1),c_4(K_2),1\}\}.
\]
\end{prop}
\begin{proof}
Let $\{\mu_1,\dotsc, \mu_{r_1}\}$ and $\{\rho_1, \dotsc, \rho_{r_2}\}$ 
be respectively systems of fundamental units for $K_1$ and $K_2$ 
as in Lemma~\ref{lem:fugood}; 
in particular we know that
\begin{equation}\label{eqn:fuheightprod}  
\prod_{j=1}^{r_1} \hrm (\mu_j) \leq c_1(K_1) R_{K_1}, \qquad\qquad
\prod_{j=1}^{r_2} \hrm (\rho_j) \leq c_1(K_2) R_{K_2}.
\end{equation}

We can write
\[
\eps_1=\zeta_1 \mu_1^{b_1}\cdots \mu_{r_1}^{b_{r_1}}, \quad
\eps_2=\zeta_2 \rho_1^{f_1}\cdots \rho_{r_2}^{f_{r_2}},
\]
where $\zeta_1$ and $\zeta_2$ are roots of unity and 
$b_1,\dotsc,b_{r_1}$, and $f_1,\dotsc,f_{r_2}$
are rational integers. Set
\[
B_1=\max\{\lvert b_1\rvert, \dotsc, \lvert b_{r_1}\rvert \}, 
\quad B_2=\max\{\lvert f_1\rvert, \dotsc,  \lvert f_{r_2}\rvert \},
\quad B=\max\{B_1, B_2,1\}.
\]
Set
$\, \alpha_0 = - \zeta_2\nu_2/(\zeta_1\nu_1)$ and $b_0=1$.
By \eqref{eqn:uniteq}, 
\[
\frac{\nu_3\eps_3}{\nu_1\eps_1}=
\alpha_0^{b_0} \mu_1^{-b_1}\cdots  
\mu_{r_1}^{-b_{r_1}}\rho_1^{f_1}\cdots \rho_{r_2}^{f_{r_2}}-1.
\]
Now choose the real or complex embedding $\sigma$ of $L$
such that $\,\lvert \sigma ((\nu_3 \eps_3)/(\nu_1 \eps_1))\rvert \,$ is minimal. 
We apply Matveev's estimate (Lemma~\ref{lem:Matveev}) 
to this \lq\lq linear form\rq\rq, obtaining 
\[
\log \left\lvert \sigma\left(\frac{\nu_3\eps_3}{\nu_1\eps_1} \right) \right\rvert >
-C(L,n) (1+\log(nB)) 
\hrm_{L,\sigma} (\alpha_0)
 \prod_{j=1}^{r_1} \hrm_{L,\sigma}(\mu_j) \prod_{j=1}^{r_2} \hrm_{L,\sigma}(\rho_j),
\]
where $n=r_1+r_2+1$. 
Using Lemma~\ref{lem:modified} and equation~\eqref{eqn:fuheightprod}  we obtain
\[
\prod_{j=1}^{r_1} \hrm_{L,\sigma}(\mu_j)
\leq (\partial_{L/K_1})^{r_1} \prod_{j=1}^{r_1} \hrm(\mu_j) 
\leq c_1(K_1) (\partial_{L/K_1})^{r_1} R_{K_1},
\]
and a similar estimate for $\prod_{j=1}^{r_2} \hrm_{L,\sigma}(\rho_j)$.
Moreover, again by Lemma~\ref{lem:modified} and Lemma~\ref{lem:elementary},
$\hrm_{L,\sigma}(\alpha_0)\leq 2H \partial_{L/L}$. 
Thus 
\[
\log \left\lvert \sigma\left(\frac{\nu_3\eps_3}{\nu_1\eps_1} \right) \right\rvert
>
-A_1 (1+\log(nB)).
\]
Now applying Lemma~\ref{lem:minsig}, we obtain that
\[
\hrm \left(\frac{\nu_3\eps_3}{\nu_1\eps_1} \right)
\leq \hrm \left(\frac{\nu_3}{\nu_1} \right)+A_1 (1+\log(nB))
\leq 2H+A_1 (1+\log(nB)).
\]
The proof is complete on observing, from Lemma~\ref{lem:Best}, that
\[
B \leq \max\{c_4(K_1),c_4(K_2),1)\} \max\{\hrm(\eps_1),\hrm(\eps_2),1\},
\]
and from Lemma~\ref{lem:elementary}, $\hrm(\nu_i \eps_i) \leq \hrm(\eps_i)+\hrm(\nu_i)
\leq \hrm(\eps)+H$.
\end{proof}

\section{Upper Bounds for the Size of \\ Integral Points on Hyperelliptic Curves}\label{sec:bounds}

We shall need the following standard sort of lemma.
\begin{lem}\label{lem:PdW}
Let $a$, $b$, $c$, $y$ be positive numbers and suppose that
\[
y \leq a+b \log(c+y).
\]
Then
\[
y \leq 2b \log{b}+2a+c.
\]
\end{lem}
\begin{proof}
Let $z=c+y$, so that $z \leq (a+c) +b \log{z}$. 
Now we apply case $h=1$ of Lemma 2.2 of \cite{PdW};
this gives $z \leq 2(b \log{b}+a+c)$, and the lemma follows.
\end{proof}

\begin{thm}\label{thm:bounds}
Let $\alpha$ be an algebraic integer of degree at least $3$,
and let $\kappa$ be a integer belonging to $K$.
Let $\alpha_1$, $\alpha_2$, $\alpha_3$ be distinct conjugates of $\alpha$
and $\kappa_1$, $\kappa_2$, $\kappa_3$ be the corresponding conjugates
of $\kappa$. Let
\[
K_1=\Q(\alpha_1,\alpha_2,\sqrt{\kappa_1 \kappa_2}),
\quad
K_2=\Q(\alpha_1,\alpha_3,\sqrt{\kappa_1 \kappa_3}),
\quad
K_3=\Q(\alpha_2,\alpha_3,\sqrt{\kappa_2 \kappa_3}),
\] 
and 
\[
L=\Q(\alpha_1,\alpha_2,\alpha_3,\sqrt{\kappa_1 \kappa_2},\sqrt{\kappa_1 \kappa_3}).
\]
Let $R$ be an upper bound for the regulators of $K_1$, $K_2$ and $K_3$.
Let $r$  be the maximum of the unit ranks of $K_1$, $K_2$, $K_3$.
Let
\[
c_j^*=\max_{1 \leq i \leq 3} c_j(K_i).
\]
Let
\[
N= 
\max_{1 \leq i,j \leq 3} 
\left\lvert\norm_{\Q(\alpha_i,\alpha_j)/\Q}(\kappa_i(\alpha_i-\alpha_j))\right\rvert^2.
\]
Let
\[
H^*= c_5^* R+
\frac{\log{N}}{\min_{1\leq i \leq 3}[K_i:\Q]} +
\hrm(\kappa).
\]
Let 
\[
A_1^*=2 H^* \cdot C(L,2r+1) 
\cdot  (c_1^*)^2   \partial_{L/L}
\cdot
\left (\max_{1 \leq i \leq 3} \partial_{L/K_i}\right)^{2r}
\cdot R^2,
\] 
and
\[
A_2^*=2H^*+A_1^*+A_1^* \log\{(2r+1) \cdot \max\{c_4^*,1\}\}.
\]
If $x \in \Z\backslash\{0\}$ satisfies $x-\alpha=\kappa \xi^2$ for some $\xi \in K$
then 
\[
\log\lvert x \rvert \leq 8A_1^* \log(4 A_1^*)+8A_2^*+H^*+20\log{2}+13 \hrm(\kappa)+19\hrm(\alpha).
\]
\end{thm}
\begin{proof}
Conjugating the relation $x-\alpha=\kappa \xi^2$ appropriately and taking
differences we obtain
\[
\alpha_1-\alpha_2=\kappa_2 \xi_2^2-\kappa_1 \xi_1^2,
\quad
\alpha_3-\alpha_1=\kappa_1 \xi_1^2-\kappa_3 \xi_3^2,
\quad
\alpha_2-\alpha_3=\kappa_3 \xi_3^2-\kappa_2 \xi_2^2.
\]
Let
\[
\tau_1=\kappa_1 \xi_1, \qquad \tau_2=\sqrt{\kappa_1 \kappa_2} \xi_2,
\qquad \tau_3=\sqrt{\kappa_1 \kappa_3} \xi_3.
\]
Observe that
\[
\kappa_1(\alpha_1-\alpha_2)=\tau_2^2-\tau_1^2,
\quad
\kappa_1(\alpha_3-\alpha_1)=\tau_1^2-\tau_3^2,
\quad
\kappa_1(\alpha_2-\alpha_3)=\tau_3^2-\tau_2^2, 
\]
and
\[
\tau_2 \pm \tau_1 \in K_1, \quad \tau_1 \pm \tau_3 \in K_2,
\quad \tau_3 \pm \tau_2 \in \sqrt{\kappa_1/\kappa_2} K_3.
\]
We claim that each $\tau_i \pm \tau_j$ can be written in the form
$\nu \eps$ where $\eps$
is a unit in one of the $K_i$ and $\nu \in L$ is an integer satisfying
$\hrm(\nu) \leq H^*$. Let us show this for $\tau_2-\tau_3$;
the other cases are either similar or easier. Note that
$\tau_2-\tau_3=\sqrt{\kappa_1/\kappa_2} \nu^{\prime\prime}$
where $\nu^{\prime\prime}$ is an integer belonging to $K_3$.
Moreover, $\nu^{\prime\prime}$ divides
\[
\sqrt{\frac{\kappa_2}{\kappa_1}}(\tau_3-\tau_2) \cdot 
\sqrt{\frac{\kappa_2}{\kappa_1}}(\tau_3+\tau_2)
=\kappa_2 (\alpha_2-\alpha_3).
\] 
Hence $\lvert \norm_{K_3/\Q}(\nu^{\prime\prime}) \rvert \leq N$. 
By Lemma~\ref{lem:improve}, we can write
$\nu^{\prime\prime}=\nu^{\prime}\eps$ where $\eps \in K_3$
and 
\[
\hrm(\nu^\prime) \leq c_5(K_3) R + \frac{\log{N}}{[K_3:\Q]}.
\]
Now let $\nu=\sqrt{\kappa_1/\kappa_2} \nu^\prime$.
Thus $\tau_2-\tau_3=\nu \eps$ where $\hrm(\nu)\leq \hrm(\nu^\prime)+
\hrm(\kappa) \leq H^*$ proving our claim.

We apply Proposition~\ref{prop} to the unit equation
\[
(\tau_1-\tau_2)+(\tau_3-\tau_1)+(\tau_2-\tau_3)=0,
\]
which is indeed of the form $\nu_1 \eps_1+\nu_2\eps_2+\nu_3\eps_3=0$
where the $\nu_i$ and $\eps_i$ satisfy the conditions of that proposition
with $H$ replaced by $H^*$.
We obtain 
\[
\hrm\left(\frac{\tau_1-\tau_2}{\tau_1-\tau_3}\right) 
\leq A_2^*+A_1^* \log\{H^*+\max\{\hrm(\tau_2-\tau_3),
\hrm(\tau_1-\tau_2)
\}\}.
\]
Observe that
\begin{equation*}
\begin{split}
\hrm(\tau_i\pm \tau_j) & \leq
\log{2}+\hrm(\tau_i)+\hrm(\tau_j)\\
& \leq \log{2}+2 \hrm(\kappa)+2 \hrm(\xi)\\
& \leq \log{2}+3 \hrm(\kappa)+ \hrm(x-\alpha)\\
& \leq 2\log{2}+3 \hrm(\kappa)+ \hrm(\alpha) + \log\lvert x \rvert,\\
\end{split}
\end{equation*}
where we have made repeated use of Lemma~\ref{lem:elementary}.
Thus
\[
\hrm\left(\frac{\tau_1-\tau_2}{\tau_1-\tau_3}\right)
\leq A_2^*+A_1^* \log(A_3^* +\log\lvert x \rvert),
\]
where $A_3^*=H^*+2\log{2}+3 \hrm(\kappa)+ \hrm(\alpha)$.

We also apply Proposition~\ref{prop} to the unit equation
\[
(\tau_1+\tau_2)+(\tau_3-\tau_1)-(\tau_2+\tau_3)=0,
\]
to obtain precisely the same bound for 
$\hrm\left(\frac{\tau_1+\tau_2}{\tau_1-\tau_3}\right)$.
Using the identity
\[
\left(\frac{\tau_1-\tau_2}{\tau_1-\tau_3}\right)\cdot
\left(\frac{\tau_1+\tau_2}{\tau_1-\tau_3}\right)=\frac{\kappa_1(\alpha_2-\alpha_1)}{(\tau_1-\tau_3)^2},
\]
we obtain that
\[
\hrm(\tau_1-\tau_3) \leq \frac{\log{2}+\hrm(\kappa)}{2}+\hrm(\alpha)+A_2^*+A_1^* \log(A_3^* +\log\lvert x \rvert).
\]
Now
\begin{equation*}
\begin{split}
\log\lvert x \rvert & \leq \log{2}+\hrm(\alpha)+\hrm(x-\alpha_1) \\
& \leq \log{2}+\hrm(\alpha)+\hrm(\kappa)+2 \hrm(\tau_1) \qquad (\text{using $x-\alpha_1=\tau_1^2/\kappa_1$})\\
&  \leq 5\log{2}+\hrm(\alpha)+\hrm(\kappa)+2 \hrm(\tau_1+\tau_3)+2 \hrm(\tau_1-\tau_3)\\
& \leq 5\log{2}+\hrm(\alpha)+\hrm(\kappa)+2 \hrm\left(\frac{\kappa_1(\alpha_3-\alpha_1)}{\tau_1-\tau_3} \right)+2 \hrm(\tau_1-\tau_3)\\
& \leq 7\log{2}+5\hrm(\alpha)+3\hrm(\kappa)+4\hrm(\tau_1-\tau_3)\\
& \leq 9\log{2}+9\hrm(\alpha)+5\hrm(\kappa)+4A_2^*+4A_1^*\log(A_3^* +\log\lvert x \rvert).
\end{split}
\end{equation*}
The theorem follows from Lemma~\ref{lem:PdW}.
\end{proof}

\section{The Mordell--Weil Sieve I}\label{sec:mwI}

The Mordell--Weil sieve is a technique that can be 
used to show the non-existence of rational points on
a curve (for example \cite{BS}, \cite{BS3}), or to help
determine the set of rational points in conjunction
with the method of Chabauty (for example \cite{BE});
for connections to the Brauer--Manin obstruction
see, for example, \cite{FlynnBM}, \cite{PoonenBM}
or \cite{StollANT}.
In this section and the next we explain how the
Mordell--Weil sieve can be used to show that any
rational point on a curve of genus $\geq 2$ is either
a known rational point or a very large rational point. 

In this section we let $C/\Q$ be a smooth projective curve
(not necessarily hyperelliptic)
of genus $g \geq 2$ and we let $J$ be its Jacobian.
As indicated in the introduction, we assume the knowledge of
some rational point on $C$; 
henceforth let $D$ be a fixed rational 
point on $C$ (or even a fixed rational divisor of degree $1$) 
and let $\jmath$ be the corresponding Abel--Jacobi map:
\[
\jmath : C \rightarrow J, \qquad P \mapsto [P-D].
\]
Let $W$ be the image in $J$ of the known rational points on $C$.
The Mordell--Weil sieve is a strategy for obtaining
a very large and \lq smooth\rq\ positive integer $B$ such that 
\[
\jmath(C(\Q)) \subseteq W+B J(\Q).
\]
Recall that a positive integer $B$ is called $A$-smooth if
all its prime factors are $\leq A$. By saying that $B$
is smooth, we loosely mean that it is $A$-smooth with $A$
much smaller than $B$.

Let $S$ be a finite set of primes, which for now we assume to be primes
of good reduction for the curve~$C$. The basic idea is to consider
the following commutative diagram.
\[ \xymatrix{C(\Q) \ar[r]^-{\jmath} \ar[d] & J(\Q)/B J(\Q) \ar[d]^-{\alpha} \\
             \displaystyle\prod_{p \in S} C(\F_p) \ar@<4pt>[r]^-{\jmath}
               & \displaystyle\prod_{p \in S} J(\F_p)/B J(\F_p)
            }
\]
The image of $C(\Q)$ in~$J(\Q)/B J(\Q)$ must then be contained in the 
subset of $J(\Q)/B J(\Q)$ of elements that map under~$\alpha$ into the
image of the lower horizontal map. If we find that this subset equals the
image of~$W$ in~$J(\Q)/B J(\Q)$, then we have shown that
\[ \jmath(C(\Q)) \subseteq W+B J(\Q) \]
as desired. Note that, at least in principle, the required computation
is finite: each set $C(\F_p)$ is finite and can be enumerated, hence
$\jmath(C(\F_p))$ can be determined, and we assume that we know explicit
generators of~$J(\Q)$, which allows us to construct the finite set
$J(\Q)/B J(\Q)$. In practice, and in particular for the application we
have in mind here, we will need a very large value of~$B$, so this
naive approach is much too inefficient. In \cite{BS} and~\cite{BS3},
the authors 
describe how one can perform this computation in a more efficient way.

One obvious improvement is to replace the lower horizontal map in the
diagram above by a product of maps
\[ C(\Q_p) \stackrel{\jmath}{\to} G_p/B G_p \]
with suitable finite quotients $G_p$ of~$J(\Q_p)$. We have used this to
incorporate information modulo higher powers of~$p$ for small primes~$p$.
This kind of information is often called ``deep'' information, as opposed
to the ``flat'' information obtained from reduction modulo good primes.

We can always force $B$ to be divisible by any given (not too big) number.
In our application we will want $B$ to kill the rational torsion subgroup
of~$J$.

\section{The Mordell--Weil Sieve II}\label{sec:mwII}

We continue with the notation of Section~\ref{sec:mwI}. 
Let $W$ be the image in $J(\Q)$ of all the known
rational points on $C$. We assume that the strategy of
Section~\ref{sec:mwI}
is successful in yielding a large \lq smooth\rq\ integer $B$
such that any point $P \in C(\Q)$ satisfies
$\jmath(P)-w \in BJ(\Q)$ for some $w \in W$, and moreover, that $B$ kills
all the torsion of $J(\Q)$.

Let
\[
\phi:\Z^r \rightarrow J(\Q), \qquad \phi(a_1,\dots,a_r)=\sum a_i D_i,
\]
so that the image of $\phi$ is simply the free part of $J(\Q)$.
Our assumption is now that 
\[
\jmath(C(\Q)) \subset W + \phi(B\Z^n).
\]

Set $L_0=B \Z^n$. We explain a method of obtaining a (very long)
decreasing sequence of lattices in $\Z^n$:
\begin{equation}\label{eqn:dec}
B\Z^n=L_0 \supsetneq L_1 \supsetneq L_2 \supsetneq \dots \supsetneq L_k 
\end{equation}
such that 
\[
\jmath(C(\Q)) \subset W + \phi(L_j)
\]
for $j=1,\dotsc,k$.  

If $q$ is a prime of good reduction for $J$ we denote
by
\[
\phi_q: \Z^r \rightarrow J(\F_q), \qquad \phi_q(a_1,\dots,a_r)=\sum a_i \tilde{D}_i,
\]
and so $\phi_q(\bl)=\widetilde{\phi(\bl)}$.

\begin{lem}\label{lem:criterion}
Let $W$ be a finite subset of $J(\Q)$, and let $L$ be a subgroup of $\Z^r$.
Suppose that $\jmath(C(\Q)) \subset W+\phi(L)$.
Let $q$ be a prime of good reduction
for $C$ and $J$. Let $L^\prime$ be the kernel of the restriction
$\phi_q \rvert_L$.
Let $\bl_1,\ldots,\bl_m$ be representatives of the {\bf non-zero} cosets
of $L/L^\prime$ and suppose that $\tilde{w}+\phi_q(\bl_i) \notin \jmath C(\F_q)$
for all $w \in W$ and $i=1,\dotsc,m$. Then 
$\jmath(C(\Q)) \subset W+\phi(L^\prime)$.
\end{lem}
\begin{proof}
Suppose $P \in C(\Q)$. Since $j(C(\Q)) \subset W+\phi(L)$, we may
write $\jmath(P)=w+\phi(\bl)$ for some $\bl \in L$. Now
let $\bl_0=\mathbf{0}$, so that $\bl_0,\dotsc,\bl_m$ represent
{\bf all} cosets of $L/L^\prime$. Then $\bl=\bl_i+\bl^\prime$
for some $\bl^\prime \in L^\prime$ and $i=0,\dotsc,m$.
However, $\phi_q(\bl^\prime)=0$, or in other words, $\widetilde{\phi(\bl^\prime)}=0$.
Hence 
\[
\jmath(\tilde{P})=\widetilde{\jmath(P)}=\tilde{w}+\phi_q(\bl)=\tilde{w}+\phi_q(\bl_i)+\phi_q(\bl^\prime)=\tilde{w}+\phi_q(\bl_i).
\]
By hypothesis, $\tilde{w}+\phi_q(\bl_i) \notin \jmath C(\F_q)$
for $i=1,\dotsc,m$, so $i=0$ and so $\bl_i=0$.
Hence $\jmath(P)=w+\bl^\prime \in W+L^\prime$ as required.
\end{proof}

\bigskip

We obtain a very long strictly decreasing sequence of 
lattices as in \eqref{eqn:dec} by repeated application of 
Lemma~\ref{lem:criterion}. However, the conditions of
Lemma~\ref{lem:criterion}
are unlikely to be satisfied for a prime $q$ chosen at random.
Here we give criteria that we have employed in practice to choose
the primes $q$.

\begin{enumerate}
\item[(I)] $\gcd(B,\# J(\F_q))> (\#J(\F_q))^{0.6}$,
\item[(II)] $L^\prime \ne L$,
\item[(III)] $\#W \cdot (\# L/L^\prime -1) < 2q$, 
\item[(IV)] $\tilde{w}+\phi_q(\bl_i) \notin \jmath C(\F_q)$
for all $w \in W$ and $i=1,\dotsc,m$.
\end{enumerate}
The criteria I--IV are listed in the order in which we check them
in practice. Criterion IV is just the criterion of the lemma.
Criterion II ensures that $L^\prime$ is strictly smaller
than $L$, otherwise we gain no new information. Although we 
would like $L^\prime$ to be strictly smaller
than $L$, we do not want the index $L/L^\prime$ to be too large
and this is reflected in Criteria I and III. Note that the
number of checks required by Criterion IV (or the lemma) is
$\#W \cdot (\# L/L^\prime -1)$. If this number is large
then Criterion IV is likely to fail. Let us look at this
in probabilistic terms. Assume that the genus of $C$ is $2$.
Then the probability that a random element of $J(\F_q)$
lies in the image of $C(\F_q)$ is about $1/q$. 
If $N=\#W \cdot (\# L/L^\prime -1)$ then the
probability that Criterion IV is satisfied is about
$(1-q^{-1})^N$. Since $(1-q^{-1})^q \sim e^{-1}$,
we do not want $N$ to be too large in comparison
to $q$, and this explains the choice of $2q$ in Criterion III.

We still have not justified Criterion I. The computation involved
in obtaining $L^\prime$ is a little expensive. Since we need
to do this with many primes, we would like a way of picking
only primes where this computation is not wasted, and in
particular $\# L/L^\prime$ is not too large. Now at every
stage of our computations, $L$ will be some element of our 
decreasing sequence~\eqref{eqn:dec} and so contained in 
$B \Z^n$. Criterion I ensures that a \lq large chunk\rq\ of
$L$ will be in the kernel of $\phi_q : \Z^n \rightarrow J(\F_q)$
and so that $\# L/L^\prime$ is not too large. The exponent $0.6$
in Criterion I is chosen on the basis of computational experience.

\section{Lower Bounds for the Size of Rational Points}\label{sec:lb}

In this section, we suppose that the strategy of Sections~\ref{sec:mwI}
and~\ref{sec:mwII} succeeded in showing that 
$\jmath(C(\Q)) \subset W+\phi(L)$ for some lattice $L$
of huge index in $\Z^r$, where $W$ is the image of $J$
of the set of known rational points in $C$.
In this section we provide a lower bound for the size of rational points
not belonging to the set of known rational points. 

\begin{lem}\label{lem:lowerbound}
Let $W$ be a finite subset of $J(\Q)$, and let $L$ be a sublattice of $\Z^r$.
Suppose that $\jmath(C(\Q)) \subset W+\phi(L)$.
Let $\mu_1$ be a lower bound for $h-\hat{h}$ as in \eqref{eqn:hdiff}.
Let
\[
\mu_2= \max \left\{ \sqrt{\hat{h}(w)} \; : \; w \in W\right\}.
\]
Let $M$ be the height-pairing matrix for the Mordell--Weil basis
$D_1,\dotsc,D_r$ and let $\lambda_1,\dotsc,\lambda_r$ be 
its eigenvalues. Let
\[
\mu_3=\min\left\{ \sqrt{\lambda_j} \; : \; j=1,\dotsc,r\right\}.
\]
Let $m(L)$ be the Euclidean norm of the shortest non-zero vector
of $L$, and suppose that $\mu_3 m(L) \geq \mu_2$. 
Then, for any $P \in C(\Q)$, either $\jmath(P) \in W$
or
\[
h(\jmath(P)) \geq \left(\mu_3 m(L) -\mu_2 \right)^2+\mu_1.
\]
\end{lem}
Note that $m(L)$ is called the 
minimum of $L$ and can 
be computed using an algorithm of Fincke and Pohst \cite{FP}.
\begin{proof}
Suppose that $\jmath(P) \notin W$. Then $\jmath(P)=w+\phi(\bl)$
for some non-zero element $\bl \in L$. In particular, if $\lVert \cdot \rVert$
denotes Euclidean norm then $\lVert \bl \rVert \geq m(L)$. 

We can write $M=N \Lambda N^t$ where $N$ is orthogonal and $\Lambda$
is the diagonal matrix with diagonal entries $\lambda_i$. Let $\bx=\bl N$ and
write $\bx=(x_1,\dotsc,x_r)$. Then
\[
\hat{h}(\phi(\bl)) =\bl M \bl^t = \bx \Lambda \bx^t \geq \mu_3^2 \lVert \bx \rVert^2
=\mu_3^2 \lVert \bl \rVert^2 \geq \mu_3^2 m(L)^2.
\]

Now recall that $D \mapsto \sqrt{\hat{h}(D)}$ defines
a norm on $J(\Q) \otimes \R$ and so by the triangle inequality
\[
\sqrt{\hat{h}(\jmath(P))} \geq \sqrt{\hat{h}(\phi(\bl))} - \sqrt{\hat{h}(w)}
\geq \mu_3 m(L) - \mu_2.
\]
The lemma now follows from \eqref{eqn:hdiff}.
\end{proof}

\noindent{\bf Remark.} We can replace $\mu_3 m(L)$ with the 
minimum of $L$ with respect to the height pairing matrix.
This is should lead to a very slight improvement. Since in practice
our lattice $L$ has very large index, computing the minimum
of $L$ with respect to the height pairing matrix may require
the computation of the height pairing matrix to very great accuracy,
and such a computation is inconvenient. We therefore prefer to work with
the Euclidean norm on $\Z^r$.

\section{Proofs of Theorems~\ref{thm:ex} and~\ref{thm:ex2}}\label{sec:example}

The equation $Y^2-Y=X^5-X$ is transformed into 
\begin{equation}\label{eqn:modified}
C \quad : \quad 2y^2=x^5-16x+8,
\end{equation}
via the change of variables $y=4Y-2$ and $x=2X$ which preserves integrality.
We shall work the model \eqref{eqn:modified}.
Let $C$ be the smooth projective genus $2$ curve with
affine model given by \eqref{eqn:modified}, and let $J$ 
be its Jacobian. Using {\tt MAGMA} \cite{MAGMA} we
know that $J(\Q)$ is free of rank $3$ with
Mordell--Weil basis given by 
\[
D_1=(0,2)-\infty,
\qquad
D_2=(2,2)-\infty,
\qquad 
D_3=(-2,2)-\infty.
\]
The {\tt MAGMA} programs used for this step are based on 
Stoll's papers \cite{StollB1}, \cite{Stoll}, \cite{StollB2}.

Let $f=x^5-16x+8$. Let $\alpha$ be a root of $f$. We shall choose for 
coset representatives of $J(\Q)/2J(\Q)$ the linear combinations
$\sum_{i=1}^{3} n_i D_i$ with $n_i \in \{0,1\}$. Then 
\[
x-\alpha=\kappa \xi^2,
\]
where $\kappa \in \cK$ 
and $\cK$ is constructed as in Lemma~\ref{lem:cK}.
We tabulate the $\kappa$ corresponding to the $\sum_{i=1}^{3} n_i D_i$  
in Table~\ref{table:bounds}.

Next we compute the bounds for $\log{x}$ given by Theorem~\ref{thm:bounds} for each
value of $\kappa$. We implemented our bounds in {\tt MAGMA}. 
Here the Galois group of $f$ is $S_5$ which implies that 
the fields $K_1$, $K_2$, $K_3$ corresponding to a particular $\kappa$ are isomorphic.
The unit ranks of $K_i$, the bounds for their regulator as given by Lemma~\ref{lem:Landau},
and the corresponding bounds for $\log{x}$ are tabulated in Table~\ref{table:bounds}.

\begin{table}
\caption{}
\begin{tabular}{||c|c|c|c|c||}
\hline\hline
coset of  &  & unit rank  & bound $R$ for  & bound for \\
$J(\Q)/2J(\Q)$ & $\kappa$ & of $K_i$ & regulator of $K_i$ & $\log{x}$ \\ 
\hline \hline
$0$ & $1$ & $12$ & $1.8 \times 10^{26}$ & $1.0\times 10^{263}$\\
$D_1$ & $-2\alpha$ & $21$ & $6.2 \times 10^{53}$ & $7.6 \times 10^{492}$\\
$D_2$ & $4-2\alpha$ & $25$ & $1.3 \times 10^{54}$ & $2.3 \times 10^{560} $\\
$D_3$ & $-4-2\alpha$ & $21$ & $3.7 \times 10^{55}$ & $1.6 \times 10^{498}$\\
$D_1+D_2$ & $-2\alpha+\alpha^2$ & $21$ & $1.0 \times 10^{52}$ & $3.2 \times 10^{487}$ \\
$D_1+D_3$ & $2\alpha+\alpha^2$ & $25$ & $7.9 \times 10^{55}$ & $5.1 \times 10^{565}$\\
$D_2+D_3$ & $-4+\alpha^2$ & $21$ & $3.7 \times 10^{55}$ & $1.6 \times 10^{498} $\\
$D_1+D_2+D_3$ & $8\alpha-2\alpha^3$ & $25$ & $7.9 \times 10^{55}$ & $5.1 \times 10^{565}$ \\
\hline\hline
\end{tabular}
\label{table:bounds}
\end{table}

A quick search reveals $17$ rational points on $C$:
\begin{gather*}
\infty, \;  (-2,  \pm 2), \; (0 , \pm 2),\; (2 , \pm 2),\; (4 , \pm 22),\;
 (6 , \pm 62), \\ (1/2, \pm 1/8), \; 
(-15/8, \pm 697/256),\;  (60 , \pm 9859).
\end{gather*}
Let $W$ denote the image of this set in $J(\Q)$.
Applying the implementation of the Mordell--Weil sieve due to Bruin and
Stoll which is explained in Section~\ref{sec:mwI} we obtain that $\jmath(C(\Q)) \subseteq W+B J(\Q)$
where
\begin{align*}
  B &= 4449329780614748206472972686179940652515754483274306796568214048000 \\
    &= 2^8 \cdot 3^4\cdot 5^3 \cdot 7^3 \cdot 11^2 \cdot 13^2 \cdot 17^2 
        \cdot 19 \cdot 23 \cdot 29 \cdot 31^2 \cdot 
        \prod_{\begin{array}{c} \scriptstyle 37 \le p \le 149 \\[-2pt] 
                                \scriptstyle p \neq 107 \end{array}} p \,.
\end{align*}
For this computation, we used \lq\lq deep\rq\rq\ 
information modulo
$2^9$, $3^6$, $5^4$, $7^3$, $11^3$, $13^2$, $17^2$, $19^2$, and
\lq\lq flat\rq\rq\ 
information from all primes $p < 50000$
such that $\#J(\F_p)$ is 500-smooth (but keeping only information
coming from the maximal 150-smooth quotient group of~$J(\F_p)$).
Recall that an integer is called \emph{$A$-smooth} if all its prime
divisors are $\le A$.
This computation took about 7 hours on a $2$ GHz Intel Core 2 CPU.

We now apply the new extension of the Mordell--Weil sieve explained in Section~\ref{sec:mwII}.
We start with $L_0=B \Z^3$ where $B$ is as above.  
We successively apply Lemma~\ref{lem:criterion} using
all primes $q < 10^6$ which are primes of good reduction 
and satisfy criteria I--IV of Section~\ref{sec:mwII}.
There
are $78498$ primes less than $10^{6}$. Of these, we discard $2$,
$139$, $449$ as they are primes of bad reduction for $C$. 
This leaves us with $78495$ primes. Of these, Criterion I fails for
$77073$ of them, Criterion II fails for $220$ of the remaining, 
Criterion III fails for $43$ primes that survive Criteria I and II,
and Criterion IV fails for $237$ primes that survive Criteria I--III.
Altogether, only $922$ primes $q<10^{6}$ satisfy Criteria I--IV and increase
the index of $L$.

The index of the final $L$ in $\Z^3$ is approximately 
$3{.}32 \times 10^{3240}$. This part of the computation 
lasted about 37 hours on a $2.8$ GHZ Dual-Core AMD Opteron.

Let $\mu_1$, $\mu_2$, $\mu_3$ be as in the notation of Lemma~\ref{lem:lowerbound}.
Using {\tt MAGMA} we find $\mu_1=2.677$, $\mu_2=2.612$ and $\mu_3=0.378$
(to $3$ decimal places). The shortest vector of the final lattice $L$
is of Euclidean length approximately $1.156 \times 10^{1080}$ 
(it should be no surprise
that this is roughly the cube root of the index of $L$ in $\Z^3$).
By Lemma~\ref{lem:lowerbound} if $P \in C(\Q)$ is not one of the $17$
known rational points then
\[
h(\jmath(P)) \geq 1.9 \times 10^{2159}.
\]
If $P$ is an integral point, then
$h(\jmath(P))=\log{2}+2\log{x(P)}$.
Thus
\[
\log{x(P)} \geq 0.95 \times 10^{2159}.
\]
This contradicts the bounds for $\log{x}$ in Table~\ref{table:bounds}
and shows that the integral point $P$ must be one of the $17$ known
rational points. This completes the proof of Theorem~\ref{thm:ex}.
The proof of Theorem~\ref{thm:ex2} is similar and we omit
the details.

The reader can find the {\tt MAGMA} programs for verifying the 
above computations at: 
{\tt http://www.warwick.ac.uk/staff/S.Siksek/progs/intpoint/}

\end{document}